\documentclass[namedreferences]{kluwer}
\usepackage{amsmath}

\newcommand\xlm{{{\mathcal X}_{\ell,m}}}
\newcommand\Z{{\mathbb Z}}
\newcommand\C{{\mathbb C}}
\newcommand\R{{\mathbb R}}
\newcommand\N{{\mathbb N}}
\newtheorem{thm}{Theorem}
\newtheorem{lem}{Lemma}
\newtheorem{cor}[thm]{Corollary}
\newtheorem{prop}[thm]{Proposition}
\newtheorem{defn}[thm]{Definition}

\begin{document}
\begin{article}
\begin{opening}
\title{Salem numbers and growth series of some hyperbolic graphs}
\date{\today}
\author{Laurent
  \surname{Bartholdi}\email{laurent@math.berkeley.edu}\thanks{The
    authors acknowledge support from the Swiss National Science
    Foundation.}}
\institute{Section de Math\'ematiques\\
  Universit\'e de Gen\`eve\\ CP 240, 1211 Gen\`eve 24\\
  Switzerland}
\author{Tullio
  G. \surname{Ceccherini-Silberstein}\email{tceccher@mat.uniroma1.it}}
\institute{Facolt\`a di Ingegneria\\
  Universit\`a del Sannio\\
  Palazzo dell'Aquila Bosco Lucarelli, Corso Garibaldi\\
  82100 Benevento,
  Italy}
\keywords{Growth series, Hyperbolic graphs and groups, Salem
  polynomials, Salem numbers}
\classification{Mathematics Subject Classification (1991)}{\parbox[t]{0.6\textwidth}{%
    \textbf{20F32} (Geometric group theory),\\
    \textbf{11R06} (PV-numbers and generalizations).}}
\begin{abstract}
  Extending the analogous result of Cannon and Wagreich for the
  fundamental groups of surfaces, we show that, for the $\ell$-regular
  graphs $\xlm$ associated to regular tessellations of the hyperbolic
  plane by $m$-gons, the denominators of the growth series (which are
  rational and were computed by Floyd and
  Plotnick~\cite{floyd-p:semiregular}) are reciprocal Salem
  polynomials.  As a consequence, the growth rates of these graphs are
  Salem numbers. We also prove that these denominators are essentially
  irreducible (they have a factor of $X+1$ when $m\equiv2\mod 4$; and
  when $\ell=3$ and $m\equiv4\mod 12$, for instance, they have a
  factor of $X^2-X+1$). We then derive some regularity properties for
  the coefficients $f_n$ of the growth series: they satisfy
  $$K\lambda^n-R<f_n<K\lambda^n+R$$
  for some constants $K,R>0$, $\lambda>1$.
\end{abstract}
\end{opening}

\section{Introduction}
We consider the graphs $\xlm$ ($\ell,m\ge3$) defined by Floyd
and Plotnick in~\cite{floyd-p:semiregular}.  These graphs are
$\ell$-regular and are the $1$-skeletons of tessellations of the sphere
(if $(\ell-2)(m-2)<4$), the Euclidean plane (if $(\ell-2)(m-2)=4$) or
the hyperbolic plane (if $(\ell-2)(m-2)>4$) by regular $m$-gons. These
tessellations were studied by Coxeter~\cite{coxeter:honey}.  When
$m=\ell=4g$ with $g$ at least two, $\xlm$ is the Cayley graph of
the fundamental group $J_g = \pi_1(\Sigma_g)$ of an orientable closed
surface $\Sigma_g$ of genus $g$, with respect to the usual set of
generators $S_g=\{a_1,b_1,\dots,a_g,b_g\}$:
$$J_g = \left< a_1, b_1, \dots, a_g, b_g \left| \prod_{i=1}^g [a_i,b_i] 
\right.\right>.$$

The growth series for $J_g$ with respect to $S_g$, namely
$$F_g(X) = \sum_{s\in J_g} X^{|s|} = \sum_{n=0}^\infty f_n X^n,$$
where $|s|=\min \{t: {s=s_1\dots s_t, s_i\in S_g\cup S_g^{-1}}\}$
denotes the word length of $s$ with respect to $S_g$ and $f_n = |\{s
\in J_g: {|s|=n}\}|$, was computed by Cannon and Wagreich
in~\cite{cannon:surfgp} and~\cite{cannon-w:surfgp} and shown to be
rational, indeed
$$F_g(X) = \frac{1 + 2X + \dots + 2X^{2g-1} + X^{2g}}{1 - (4g-2)X -
  \dots - (4g-2)X^{2g-1} + X^{2g}};$$
moreover they showed that the
denominator is a Salem polynomial. It was later shown by
Floyd~\cite{floyd:salem} and Parry~\cite{parry:salem} that the
denominators of the growth series of Coxeter groups are also Salem
polynomials.

In~\cite{floyd-p:fuchsian} and~\cite{floyd-p:semiregular}, Floyd and
Plotnick, among other things, extended the calculations of Cannon and
Wagreich to the family $\xlm$.  Fixing arbitrarily a base point
$*\in V(\xlm)$ and denoting by $|x|$ the graph distance between
the vertices $x$ and $*$, they obtained the following formul\ae\ for the
growth series $F_{\ell,m}(X) = \sum_{x\in V(\xlm)} X^{|x|}$; for
$m$ even, say $m = 2w$:
\begin{equation}
  F_{\ell,m}(X) = \frac{1 + 2X + \dots + 2X^{w-1} + X^w}{
  1 - (\ell-2)X - \dots - (\ell-2)X^{w-1} + X^w}\label{ieq:even}
\end{equation}
and, for $m$ odd, say $m = 2w+1$:
\begin{equation}
  F_{\ell,m}(X) = \frac{1 + 2X + \dots + 2X^{w-1} + 4X^w + 2X^{w+1} + 
    \dots + 2X^{2w-1} + X^{2w}}{
    1 - (\ell-2)X - \dots - (\ell-4)X^w - \dots - (\ell-2)X^{2w-1} + X^{2w}}.\label{ieq:odd}
\end{equation}


Our main result is the following:
\begin{thm}\label{thm:main}
  The denominators of the growth series $F_{\ell,m}$ are reciprocal
  Salem polynomials. After simplification by $X+1$ in case $m\equiv 2
  \mod 4$, they are irreducible, except for the following exceptional
  cases:
  \begin{enumerate}
  \item $\ell=3$ and $m\equiv 4\mod 12$; there is a factor of
    $X^2-X+1$;
  \item $\ell=3$ and $m\equiv 3\mod 6$; there is a factor of
    $X^2+X+1$;
  \item $\ell=3$ and $m\equiv 5\mod 20$; there is a factor of
    $X^4-X^3+X^2-X+1$;
  \item $\ell=4$ and $m\equiv 3\mod 8$; there is a factor of
    $X^2+1$;
  \item $\ell=5$ and $m\equiv 3\mod 12$; there is a factor of
    $X^2-X+1$.
  \end{enumerate}
\end{thm}
This theorem extends the results in Section~4
of~\cite{floyd-p:semiregular} which show, for even $m$, that the
denominator of $F_{\ell,m}$ is a product of an irreducible Salem
polynomial and distinct cyclotomic polynomials.

\begin{cor}
  The growth rates of the graphs $\xlm$ are Salem numbers.
\end{cor}

We thus obtain more precise information about the growth coefficients:
\begin{cor}\label{cor:main}
  If $F_{\ell,m}(X) = \sum_{n\ge0}f_n X^n$, then there exist constants
  $K>0$, $\lambda>1$ and $R>0$ such that
  $$K\lambda^n-R<f_n<K\lambda^n+R$$
  holds for all $n$. Moreover $\lambda$ is a Salem number.
\end{cor}

This improves, for the graphs $\xlm$, on a general result by
Coornaert~\cite{coornaert:mesures} for non-elementary hyperbolic
groups, which asserts that there exist constants $\lambda>1$ and
$0<K_1<K_2$ such that
$$K_1\lambda^n < f_n < K_2\lambda^n$$
for all $n$.

Note that Corollary~\ref{cor:main} does not hold for a general
presentation of a hyperbolic group, nor even of a virtually free
group. Consider for instance the modular group $PSL_2(\Z)=\langle
a,b|\,a^2=b^3=1\rangle$.  The growth coefficients $f_n$ satisfy
$$K\phi^n-2<f_n<K\phi^n+2\qquad\text{with
  }K=\frac{3+\sqrt5}{\sqrt5},\phi=\frac{1+\sqrt5}{2}$$
for the
generating set $\{a,ab,(ab)^{-1}\}$ (note that $\phi$ is a Salem
number), but for the generating set $\{a,b,b^{-1}\}$ are
$$f_n = \begin{cases}2\cdot\sqrt2^n & \text{ if $n$ is even,}\\
  3/\sqrt2\cdot\sqrt2^n & \text{ if $n$ is odd.}\end{cases}$$
These computations are due to Mach\`\i; see~\cite[VI.7]{harpe:cgt}. It
follows that Corollary~\ref{cor:main} does not extend to arbitrary
hyperbolic group presentations (indeed, $\sqrt 2$ is not even a Salem
number).

\section{Salem Polynomials}
We recall a few facts on Salem polynomials; one might also
consult~\cite[\S~5.2]{bertin-:salem} or the original
paper~\cite{salem:ps}.  A polynomial $f(X) = f_0 + f_1X + \dots + f_n
X^n$ ($f_n\neq0$) is \emph{reciprocal} if $f_i = f_{n-i}$ for all
$i=0,1,\ldots,n$; equivalently if $X^nf(X^{-1}) = f(X)$.

A polynomial $f(X)\in\Z[X]$ is a \emph{Salem polynomial} if it is
monic, admits exactly one root $\lambda$ of modulus $|\lambda|>1$, and
this root is simple; this $\lambda$ is then necessarily real. If
moreover $f$ is reciprocal, then $1/\lambda$ is also a root of $f$ and
these are the only roots of $f$ not on the unit circle.
   
If $f$ is a reciprocal polynomial of odd degree, then $f(-1)=0$, so
$X+1$ divides $f$. There is therefore no real limitation in
considering reciprocal Salem polynomials of even degree.

If $f$ is a Salem polynomial and $g$ is a cyclotomic polynomial, then
$fg$ is again a Salem polynomial, and reciprocally any Salem
polynomial can be factored as a product of an irreducible Salem
polynomial and cyclotomic polynomials.

A real number $\lambda$ is called a \emph{Salem number} if
$\lambda>1$, $\lambda$ is an algebraic integer, and all its conjugates
except $\lambda^{\pm1}$ have absolute value~$1$; equivalently if
$\lambda>1$ is the root of a Salem polynomial.

For instance, the reciprocal Salem polynomials of degree $2$ are the
$X^2-aX+1$ for all $a\in\Z$ with $a\ge3$. The corresponding Salem
numbers are the $(a+\sqrt{a^2-4})/2$. The Salem polynomials of degree
$2$ are all the $X^2-aX+b$ subject to $a>|1+b|$ and $b\neq0$.

Denoting by
$$\Phi(z) = \frac{z-i}{z+i}$$
the Cayley transform~\cite[13.17]{rudin:fa} mapping the
extended real axis $\R\cup\{\infty\}$ onto the unit circle
$\mathbb T$, we give the following
\begin{defn}
  Let $f$ be a polynomial of degree $n$. Its \emph{Cayley transform}
  is the polynomial
  $$C(f)(X) = (X+i)^n f(\Phi(X)).$$
\end{defn}

Note that if $f$ is real and reciprocal, its Cayley transform will
again be real.
 
The proof of the following characterization of reciprocal Salem
polynomials is straightforward:
\begin{thm}
  Let $f$ be a monic integral reciprocal polynomial of degree $n$.
  Then $f$ is a Salem polynomial if and only if the polynomial $C(f)$
  has exactly $n-2$ real roots (its last two roots are then complex
  conjugate).
\end{thm}

\section{The Denominators of the Growth Series $F_{\ell,m}$}\label{sec:denom}
The objective of this section is to prove that the denominator of
$F_{\ell,m}$ is (after a possible division by $X+1$ depending on the
parity of its degree $m$) a reciprocal Salem polynomial.  For this
purpose define the following auxiliary polynomials, for $a\in\Z$ and
$b,k\in\N$:
\begin{align}
  p_b(X) &= 1 + 2X + 2X^2 + \dots + 2X^{b-1} + X^b = (1-X^b)\frac{1+X}{1-X},\\
  q_b(X) &= 1 + X^b,\\
  r_{a,b}(X) &= 1 + aX + aX^2 + \dots + aX^{b-1} + X^b,\label{ieq:rab}\\
  r_{a,b;\,k}(X) &= 1 + aX + aX^2 + \dots + aX^{b-1} + X^b + kX^{b/2}\quad\text{($b$ even)}.
\end{align}

\begin{lem}\label{lem:rolle}
  Let $f,g:\R\to\R$ be two continuous functions, such that the
  following holds: $f$ has $t$ real zeroes $\rho_1<\dots<\rho_t$, and
  $g$ has a transverse zero in each interval $]\rho_i,\rho_{i+1}[$ and
  no other zero in $[\rho_1,\rho_t]$. Then for any
  $\alpha,\beta\in\R\setminus\{0\}$ the function $\alpha f+\beta g$
  has at least $t-1$ real zeroes.
\end{lem}
\begin{pf}
  The function $\alpha f+\beta g$ has opposite signs at $\rho_i$ and
  $\rho_{i+1}$, so has a zero in $]\rho_i,\rho_{i+1}[$.
\end{pf}

\begin{lem}\label{lem:rabsalem}
  If $a\in\Z$, $b\in2\N$, and $2 + a(b-1) < 0$, then $r_{a,b}$ is a
  reciprocal Salem polynomial.
\end{lem}
\begin{pf}
  First write
  $$r_{a,b}(X) = \frac a2p_b(X) + (1-\frac a2)q_b(X).$$
  As $p_b$ and $q_b$ are symmetric polynomials, their Cayley
  transforms are polynomials with real coefficients.  The zeroes of
  $p_b$ are $-1$ and the $e^{2i\pi\frac kb}, k\in\{1,\dots,b-1\}$; the
  zeroes of $q_b$ are the $e^{i\pi\frac{2k+1}{b}},
  k\in\{0,\dots,b-1\}$.  As these zeroes are interleaved on the unit
  circle, the zeroes of their Cayley transforms will be interleaved on
  the real axis, and we may apply Lemma~\ref{lem:rolle} to conclude
  that $r_{a,b}(X)$ has at least $b-1$ zeroes on the unit circle. It
  has precisely $b-1$ zeroes there because $r_{a,b}(1) = 2+a(b-1) < 0$
  and $\lim_{x \rightarrow \infty} r_{a,b}(x)=+\infty$, so that
  $r_{a,b}$ has a zero in $]1,\infty[$.
\end{pf}

\begin{lem}\label{lem:perturb}
  Let $f(X)\in\Z[X]$ be a reciprocal Salem polynomial.  Let
  $g(X)\in\Z[X]$ be a polynomial of degree less than $\deg(f)$, such
  that $f+g$ is reciprocal. Consider for $\epsilon\in\R$ the
  perturbation $f_\epsilon = f + \epsilon g\in\R[X]$.  Let $k\in\N$ be
  such that $f_\epsilon$ has only simple roots for all
  $\epsilon\in[0,k]$.  Then $f_k$ is a reciprocal Salem polynomial.
\end{lem}
\begin{pf}  
  Let $F$ and $F_\epsilon$ be the Cayley transforms of $f$ and
  $f_\epsilon$.  Then $F$ has real roots except two which are complex
  conjugate, and has real coefficients. By assumption, the
  discriminant of $F_\epsilon$ does not change its sign on $[0,k]$ and
  thus $F_k$ has real roots except two which are still complex
  conjugate. Taking the inverse Cayley transform yields $f_k$ which
  has all its roots on the unit circle except two, and thus is a Salem
  polynomial.
\end{pf}

\begin{lem}\label{lem:rabksalem}
  If $a\in\Z$, $b\in2\N$ with $2+a(b-1)<0$, and $k\in\N$ with
  $k\le\min\{2-a,-2-a(b-1)\}$, then $r_{a,b;\,k}$ is a reciprocal
  Salem polynomial.
\end{lem}
\begin{pf}
  Denote $\min\{2-a,-2-a(b-1)\}$ by $K$.  In view of
  Lemma~\ref{lem:perturb}, it suffices to prove that
  $r_{a,b;\;\epsilon}=r_{a,b}+\epsilon X^{b/2}$ is simple for all
  $\epsilon\in[0,K[$. For this purpose consider the function
  $f(X)=X^{-\frac b2}r_{a,b}(X)=X^{-\frac
  b2}\left(1+X^b+a\frac{1-X^b}{1-X}\right),$ considered as a function
  on the circle ($f:\mathbb T\to\C$).
  
  Since $f$ is real and $f(X^{-1})=f(X)$, it satisfies $f(\mathbb
  T)\subset\R$: for any $\xi\in\mathbb T$,
  $$\overline{f(\xi)}=\overline f(\overline\xi)=f(\xi^{-1})=f(\xi).$$
  We show that $r_{a,b;\;\epsilon}$ is simple by showing that
  $f+\epsilon$ has only simple zeroes on $\mathbb T$, or equivalently
  that $f$ attains values greater than $\epsilon$ between its zeroes
  on $\mathbb T$.
  
  Since $r_{a,b}$ is a reciprocal Salem polynomial, it has $b-2$
  zeroes on $\mathbb T$, so $f$ has also $b-2$ zeroes on $\mathbb
  T$. Consider the points $\xi_j=e^{2i\pi j/b}\in\mathbb T$, for
  $j\in\{1,\dots,b-1\}$. We have
  $$f(\xi_j)=\xi_j^{-\frac
    b2}\left(\xi_j^b+1+a\frac{\xi_j^b-\xi_j}{\xi_j-1}\right)=(-1)^j(2-a),$$
  so the zeroes of $f$ are separated by extrema of at least
  $\pm(2-a)$.

  Finally, since $f(1)+\epsilon\le \epsilon-K<0$, the two real zeroes
  of $f+\epsilon$ remain always separated by the unit circle. The
  zeroes of $f+\epsilon$ (and thus of $r_{a,b;\;\epsilon}$) on
  $\mathbb T$ therefore remain simple.
\end{pf}

We now turn to the factorization of the $r_{a,b}$. Let $Q$ be any
reciprocal Salem polynomial, let $\lambda$ be the Salem number
associated to $Q$, and let $S$ be the minimal polynomial of $\lambda$.
Then we have a factorization $Q=ST$, where $T$, having only roots on
the unit circle, is a product of cyclotomic polynomials, in virtue of
the theorem of Kronecker~\cite[Vol.\ III, Part~I, pages
47--110]{kronecker:werke}. We show that this show that this factor $T$
is either $1$ or $X+1$, depending on the parity of $b$:
\begin{prop}\label{prop:irred}
  Suppose $|a-1|\ge2$. Then the only cyclotomic polynomials dividing
  the $r_{a,b}$ defined in~(\ref{eq:rab}) are $X+1$ when $b$ is odd,
  and $X^2-X+1$ when $a=-1$ and $b\equiv2\mod 6$, with the exception
  $a=-2,b=2$ when $r_{a,b}=(X-1)^2$ is not a Salem polynomial.
\end{prop}
\begin{pf}
  Set $f(X)=r_{a,b}(X)$. Clearly $f(-1)=0$ precisely when $b$ is odd,
  and $f(1)\neq0$. Suppose now that $\eta$ is a root of unity of order
  $n>2$ satisfying $f(\eta)=0$. We may suppose, by direct computation,
  that $n$ does not divide $b-1$. Let $\xi$ be an algebraic conjugate
  of $\eta$ satisfying $|\xi^{b-1}-1|\ge1$. Since
  $$\left|\frac{\xi^{b+1}-1}{\xi^{b-1}-1}\right|=\left|\xi\left(1-a+\frac{f(\xi)}{\xi+\dots+\xi^{b-1}}\right)\right|=|a-1|,$$
  and $|\xi^{b+1}-1|\le2$, we have $|a-1|\le2$. Now in case $|a-1|=2$
  and $\xi$ is such a root, we must have $|\xi^{b+1}-1|=2$ whence
  $\xi^{b+1}=-1$. Then we have $\xi^b-\xi=(1-a)/2=\pm1$, or, after
  multiplication by $\xi$ and simplification, $\xi^2\pm\xi+1=0$, so
  $\xi$ is of degree $3$ or $6$. We
  check by substitution in $f$ that indeed $\xi^2-\xi+1$ divides $f$
  when $b\equiv2\mod 6$. We have shown that the only roots of unity of
  degree greater than $2$ are the sixth, and that they occur only in
  very special cases.
\end{pf}

A similar, but more complicated, result holds for $r_{a,b;\;2}$; then
there are special cases for $-3\le a\le-1$, with factors $T$ depending
on the value of $b$ modulo $4$, $6$ and $10$, as mentioned in the
statement of the theorem; we omit the uninteresting details and quote
the result without proof.

\begin{pf}[Proof of Theorem~\ref{thm:main}]
  If $m$ is even the denominator of $F_{\ell,m}$ is $r_{2-\ell,\frac
    m2}$, by~(\ref{eq:even}), and thus is a (reciprocal) Salem
  polynomial by Lemma~\ref{lem:rabsalem}. If $m$ is odd, this
  denominator is $r_{2-\ell,m-1;\,2}$, by~(\ref{eq:odd}), and thus is
  a (reciprocal) Salem polynomial, by Lemma~\ref{lem:rabksalem}. It is
  irreducible, by Proposition~\ref{prop:irred}.
\end{pf}

\section{The Growth of the Graphs $\xlm$}
Recall the following well-known fact from~\cite[p.~341]{graham-k-p:concmath}:
\begin{lem}\label{fact:gkp}
  Let $f(X) = P(X)/Q(X) = \sum_{n\ge0} f_n X^n$ be a rational function
  of $X$, where $P$ and $Q$ are complex polynomials and $Q(X) =
  \prod_{i=1}^r(X-\alpha_i)^{\nu_i}$ with distinct $\alpha_i\in\C$,
  namely $\alpha_i\neq\alpha_j$ when $i\neq j$.  Then there exist
  polynomials $R_1, \dots, R_r\in\C[X]$ such that $f_n =
  \sum_{i=1}^r {R_i(n)} / {\alpha_i^n}$. Moreover the degree of $R_i$
  is strictly smaller than $\nu_i$, for all $i$. In particular, if all
  poles of $f$ are simple, then the $R_i$ are constant.
\end{lem}

From this, we derive the following result:
\begin{thm}
  Let $$f(X) = \frac{P(X)}{Q(X)} = \sum_{n\ge0} f_n X^n$$
  be a
  rational function of $X$, where $P\in\Z[X]$ and $Q$ is a Salem
  polynomial. Then there exist a constant $K>0$ and a polynomial $R$
  such that for all $n$ we have
  $$K\lambda^n - R(n) < f_n < K\lambda^n + R(n),$$
  where $\lambda>1$
  is the Salem number associated to $Q$.  The degree of $R$ is
  strictly less than the maximal multiplicity of $f$'s poles.  Thus if
  moreover all poles of $f$ are simple, then there exist constants
  $\lambda>1$, $K$ and $R$ such that
  $$K\lambda^n - R < f_n < K\lambda^n + R.$$
\end{thm}
\begin{pf}
  Apply Lemma~\ref{fact:gkp} to $f$ to obtain polynomials $R_1,\dots,
  R_r$.  Without loss of generality we may assume that $\alpha_r$ is
  the only pole of $f$ inside the unit circle. Thus $K := R_r$ is a
  constant.  Set also $\lambda = 1/\alpha_r$. Writing $R_i(n) = \sum
  b_{ij}n^j$, we define polynomials $S_i$ by
  $S_i(n) = \sum |b_{ij}|n^j$, and we let
  $$R(n) := \sum_{i=1}^{r-1} S_i(n).$$
  Then $|R_i(n)/\alpha_i^n| \le S_i(n)$, and
  $$|f_n - K\lambda^n| = |f_n - \frac{R_r(n)}{\alpha_r^n}| \le R(n).$$
  
  If all poles of $f$ are simple, then the $R_i$ are constants and so
  is $R$.
\end{pf}

Corollary~\ref{cor:main} follows from the previous theorem.

\acknowledgements
The authors express their thanks to Pierre de la Harpe, Alexander
Borisov and Kurt Foster for their interest and generous contribution.
\bibliographystyle{klunamed}
\bibliography{mrabbrev,people,bartholdi,grigorchuk,math}
\end{article}
\end{document}